\documentclass[11pt]{article}
\usepackage{enumerate}
\usepackage{amsthm,amsmath,amssymb}
\usepackage{graphicx}
\usepackage{lineno}
\usepackage[colorlinks=true,citecolor=black,linkcolor=black,urlcolor=blue]{hyperref}
\usepackage[english]{babel}
\usepackage{amsfonts}
\usepackage{epsfig, subfigure}
\usepackage{amscd,latexsym}
\usepackage{url}
\usepackage{color}
\usepackage{authblk}
\usepackage{subfigure}

\usepackage{multicol}  

\usepackage{enumitem}
\usepackage{algorithm}%
\usepackage{algorithmicx}%
\usepackage{algpseudocode}%

\theoremstyle{plain}
\newtheorem{theorem}{Theorem}
\newtheorem{lemma}[theorem]{Lemma}
\newtheorem{corollary}[theorem]{Corollary}
\newtheorem{proposition}[theorem]{Proposition}
\newtheorem{conj}[theorem]{Conjecture}

{\itshape}{\rmfamily}

\newtheorem{definition}[theorem]{Definition}

\newcommand{\ds}{\displaystyle }
\newcommand{\ma}[2]{ $ \left( \begin{array} {#1} #2 \end{array} \right) $ }

\newcommand{\ecc}{\rm ecc} 

\newcommand{\bit}{\begin{itemize}}
\newcommand{\eit}{\end{itemize}}
\newcommand{\ben}{\begin{enumerate}}
\newcommand{\een}{\end{enumerate}}

\begin{document}

\title{Not every  graph can be reconstructed from its boundary distance matrix}

\author[1]{Jos\'e C\'aceres\thanks{jcaceres@ual.es}}
\author[2]{Ignacio M. Pelayo\thanks{ignacio.m.pelayo@upc.edu}}

\affil[1]{Departmento de Matem\'aticas, Universidad de Almer\'{\i}a, Spain}
\affil[2]{Departament de Matem\`atiques, Universitat Polit\`ecnica de Catalunya, Spain}

\date{}

\maketitle
          
\abstract{
A vertex $v$ of a connected graph $G$ is said to be a boundary vertex of $G$ if for some other 
vertex $u$ of $G$, no neighbor of $v$ is further away from $u$  than $v$. 
The boundary  $\partial(G)$ of  $G$ is the set of all of its boundary vertices.

The boundary distance matrix $\hat{D}_G$ of a graph $G=([n],E)$ is the square matrix of order $\kappa$, being $\kappa$ the order of $\partial(G)$, such that for every $i,j\in \partial(G)$,  $[\hat{D}_G]_{ij}=d_G(i,j)$.

In a recent paper [doi.org/10.7151/dmgt.2567], it was shown that if a graph $G$ is either a block graph or a unicyclic graph, then $G$ is uniquely determined by the boundary distance matrix $\hat{D}_{G}$  of $G$, and it was also conjectured that this statement holds for every connected graph $G$, whenever both the order $n$ and the boundary (and thus also the boundary distance matrix) of $G$ are prefixed.

After proving  that this conjecture is  true for several graph families, such as being of diameter 2, having  order at most $n=6$ or being Ptolemaic, we show that this statement does not hold when considering, for example, either  the family of split graphs of diameter 3 and order at least $n=10$ or the family of distance-hereditary graphs of order at least $n=8$.
}

\vspace{+.1cm}\noindent \textbf{Keywords:} Metric dimension.

\vspace{+.1cm}\noindent \textbf{AMS subject classification:} 05C50, 15A18, 05C69.



\section{Introduction }

Consider a connected, simple, and undirected graph $G$. 
The boundary of $G$, denoted by $\partial(G)$, is defined as the subset of vertices whose eccentricity is maximal among their neighbors. 
This concept was first introduced by Chartrand et al. in their foundational work~\cite{cejz03}, and since then, it has found numerous applications in graph theory, particularly in the study of convexity and metric localization in graphs.

In a recent study~\cite{cp24}, the authors proposed a conjecture stating that the distance matrix of the boundary vertices, combined with the total number of vertices in the graph, is sufficient to reconstruct the entire graph. 
This is not the first time that a similar question is posed by far. Various approaches include reconstructing
metric graphs from density functions~\cite{dww18}, road networks from a set of
trajectories~\cite{aw12}, graphs utilizing shortest paths or distance oracles~\cite{kmz18}, labeled
graphs from all r-neighborhoods~\cite{mr17}, or reconstructing phylogenetic trees~\cite{bc09} has been proposed previously. Of particular note is the graph reconstruction conjecture~\cite{k42,u60} which states the possibility of reconstructing any graph on at least three vertices (up to isomorphism) from the multiset of all unlabeled subgraphs obtained through the removal
of a single vertex.

Going back to the conjecture proposed in~\cite{cp24}, it holds true for trees, since the boundary coincides with the set of leaves. Additionally, the same work establishes the validity of the conjecture for block graphs and unicyclic graphs.

However, in this paper, we demonstrate that, unfortunately, the conjecture does not hold for general graphs. 
Despite this, the question of determining for which graph classes the conjecture remains valid is of significant theoretical and computational interest. 
We prove that for two well-known families —Ptolemaic graphs and 2-diameter graphs— the distance matrix of the boundary vertices, along with the vertex count, indeed provides enough information to reconstruct the original graph. 
On the other
 hand, we show that this property fails for  bipartite graphs, distance-hereditary graphs and split graphs of diameter 3. The case of interval graphs presents a particularly intriguing challenge: although we are unable to definitively conclude whether the conjecture holds for this class, we provide a structural characterization of the boundary vertices, which may serve as a stepping stone for future investigations.

\vspace{.3cm}
All the graphs considered are undirected, simple, finite, and (unless otherwise stated) connected.
If $G=(V,E)$ is a graph of order $n$ and size $m$, it means that $|V|=n$ and $|E|=m$. 
Unless otherwise specified, $V=[n]=\{1,\ldots,n\}$.

\vspace{.3cm}
Let $v$ be a vertex of a graph $G$.
The \emph{open neighborhood} of $v$ is $\displaystyle N(v)=\{w \in V(G) :vw \in E\}$, and the \emph{closed neighborhood} of $v$ is $N[v]=N(v)\cup \{v\}$ .
The \emph{degree} of $v$ is $\deg(v)=|N(v)|$.
The minimum degree  (resp., maximum degree) of $G$ is $\delta(G)=\min\{\deg(u):u \in V(G)\}$ (resp., $\Delta(G)=\max\{\deg(u):u \in V(G)\}$).
If $\deg(v)=1$, then $v$ is said to be a  \emph{leaf}.
Given $u,v\in V(G)$, they are called either \emph{true twin} vertices or \emph{false twin} vertices if either  $N[u]=N[v]$ or $N(u)=N(v)$.

\vspace{.3cm}
Let $W\subseteq V(G)$ be a subset of vertices of  $G$.
The  \emph{closed neighborhood} of $W$ is $N[W]=\cup_{v\in W} N[v]$.
The subgraph of $G$ induced by $W$, denoted by $G[W]$, has $W$ as vertex set and $E(G[W]) = \{vw \in E(G) : v \in W,w \in W\}$.
If $G[W]$ is a complete graph, then it is said to be a \emph{clique}\index{clique} of $G$.

\vspace{.3cm}
Given a pair of vertices $u,v$ of a graph $G$, a $u-v$ \emph{geodesic} lies on a  $u-v$ shortest  path, i.e., a  path joining $u$ and $v$ of minimum order. 
Clearly,  all $u-v$ geodesics have the same length, and it is called the \emph{distance} between vertices $u$ and $v$ in $G$, denoted by $d_G(u,v)$, or simply by $d(u,v)$, when the context is clear. 
The \emph{eccentricity} $\ecc({\it v})$  of a vertex $v$ is the distance to a farthest vertex from $v$. 
The \emph{radius} and \emph{diameter} of $G$ are respectively, the minimum and maximum eccentricity of its vertices and are denoted as ${\rm rad}(G)$ and ${\rm diam}(G)$. 
A vertex $u\in V(G)$  is a \emph{central} vertex of $G$ if $\ecc({\it u})={\rm rad}(G)$, and it is called  a  \emph{peripheral} vertex of $G$  if $\ecc({\it u})={\rm diam}(G)$.
The  \emph{distance matrix} $D_G$ of a graph $G=(V,E)$ with $n$ vertices is the square matrix of order $n$ such that,  for every $i,j \in [n]$,  $[D_G]_{ij}=d(i,j)$.

\newpage
\vspace{.2cm}
Let $S$ be a subset of vertices of order $k$ of a graph $G=([n],E)$. 
It is denoted by $D_{S,V}$ the submatrix of $D_G$ of order $k\times n$ such that for every $i\in S$ and for every $j\in V$, $[D_{S,V}]_{ij}=d(i,j)$.
Similarly, the so-called \emph{$S$-distance matrix} of $G$, denoted by $D_{S}$, is the square submatrix of $D_G$ of order $k$ such that for every $i,j\in S$, 
$[D_{S}]_{ij}=d(i,j)$.

\vspace{.2cm}
For additional details and information on basic graph theory, we refer the reader to \cite{clz16}. 
The rest of this paper is organized as follows: in the next section it is established the conjecture with the definitions and facts that are related to it. Two subsections are devoted to Ptolemaic graphs and interval graphs. 
In Section~\ref{nbdm}, we provide various examples of graphs for which the conjecture is not true, providing the limits of the conjecture. 

\section{BDM graphs }\label{bdm}

A vertex $v$ of a graph $G$ is said to be a \emph{boundary vertex} 
of a vertex $u$ if no neighbor of $v$ is further away from $u$ than $v$, i.e., if for every vertex $w\in N(v)$, 
$d(u,w) \le d(u,v)$. 
The set of boundary vertices of a vertex $u$ is denoted by $\partial(u)$.
Given a  pair of vertices $u,v\in V(G)$ if $v\in \partial(u)$, then $v$ is also said to be \emph{maximally distant} from $u$.
The \emph{boundary} of $G$, denoted by $\partial(G)$, is the set of all of its boundary vertices, i.e., 
$\ds \partial(G)=\cup_{u \in V(G)}\partial(u)$ \cite{cejz03}.
The \emph{boundary number} of $G$, denoted by $\kappa(G)$, is the cardinality of all of its boundary, i.e.,
$\kappa(G)=|\partial(G)|$.
If $S=\partial(G)$, then the  \emph{$S$-distance matrix}  $D_{S}$ is also denoted by $\hat{D}_G$, being  called the \emph{boundary distance matrix} of $G$.

The boundary of a graph exhibits a number of interesting properties, like being geodetic~\cite{chmpps06}, a resolving set~\cite{hmps13} and also a strong resolving set.

\begin{theorem}[\cite{ryko14}]
The boundary $\partial(G)$ of every graph $G$ is a strong resolving set of $G$.
\label{bsds}
\end{theorem}

The  concept of strong resolving set was first defined by Seb\H o and Tannier~\cite{st04} in 2003, and later studied in~\cite{op07}. 
They were interested in extending isometric embeddings of subgraphs into the whole graph and, to ensure that, they defined a \emph{strong resolving set} of a graph $G$ as a subset $S\subseteq V(G)$ such that for any pair $x,y\in V(G)$ there is an element $v\in S$ such that  there exists either a $x-v$ geodesic that contains $y$ or a $y-v$ geodesic containing $x$. 
As a consequence of the definition, it was implicitly  mentioned in some papers \cite{k20,ryko14,st04} and proved in \cite{cp23} that it only suffices to know the distances from the vertices of a strong resolving set to the rest of  nodes, to uniquely determine the graph.

\begin{theorem}[\cite{cp23}]
Let $S$  be a proper subset of vertices  of a graph $G=(V,E)$.
Then, $S$ is a strong resolving set if and only if
$G$ is uniquely determined by the distance matrix $D_{S,V}$.
\label{sdim.dmatrix}
\end{theorem}

Hence, as a direct consequence of these facts, the following result holds.

\begin{corollary}
Let $G=(V,E)$ be a graph.
Let $\partial(G)$ the set of boundary vertices of $G$.
Then, $G$ is uniquely determined by the distance matrix $D_{\partial(G),V}$.
\end{corollary}

Starting from this property, it was conjectured  in \cite{cp24} that every graph $G$ of order $n$ and boundary number 
$\kappa$  could also be uniquely determined by the boundary distance matrix $\hat{D}_G$, provided that both the order $n$ and the boundary   $\partial(G)$ be prefixed.
In the same paper, it was proved this claim to be true  for trees, for unicyclic graphs  and also for block graphs.
In the next section, it is shown that there are graphs not satisfying this property.
For this reason, we introduce the following definition.

\begin{definition}
A graph $G$ of order $n$ and boundary number 
$\kappa$ is called {\bf boundary-distance-matrix constructible}, also {\bf BDM-constructible}, or simply a {\bf BDM} graph, 
if it is uniquely determined by the distance matrix $\hat{D}_{G}$ of the boundary $\partial(G)$ of $G$, provided that both  $n$ and $\partial(G)$ be prefixed.
\label{bdmdef}
\end{definition}

\begin{figure}[ht]
\begin{center}
\includegraphics[width=0.5\textwidth]{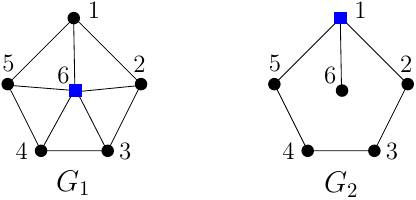}
\caption{
A pair of graphs of order  $n=6$,  boundary number $\kappa=5$, with the same $[5]$-distance matrix.}
\label{n6k5}
\end{center}
\end{figure}

At this point, it is important to notice that in the previous definition, it is not enough to prefix both the order $n$ and the boundary number $\kappa$.
For example, for both of the graphs, $G_1$ and $G_2$, displayed in Figure \ref{n6k5},  $n=6$ and  $\kappa=5$, they have the same  $[5]$-distance matrix  $D_{\stackrel{}{[5]}}$, but $\partial(G_1)=[5]$ and $\partial(G_2)=\{1,2,3,4,6\}$.

Next, we show that all graphs of diameter 2  are BDM graphs.

\begin{lemma} 
\label{b=n} 
Let $G$ be a graph  of order $n$ and boundary number $\kappa$.

\begin{enumerate}[label=\rm \bf(\arabic*)]

\item
If $v\in V(G)$ and ${\rm ecc}(v)={\rm diam}(G)$, then $v\in \partial(G)$.

\item
If ${\rm rad}(G)={\rm diam}(G)$, then $\kappa= n$.

\item
If ${\rm diam}(G)=2$, then $n-1 \le \kappa \le n$.
\newline
Moreover, $\kappa=n-1$ if and only if $G$ contains a unique central vertex.
\end{enumerate}
\begin{proof}
\begin{enumerate}[label=\rm \bf(\arabic*)]
\item
If  $w\in V(G)$ is an antipodal  vertex of $v$, i.e., if   $d(w,v)={\rm diam}(G)$, then $v\in \partial(w)$.
\item
If every vertex of $G$ is peripheral, then according to the previous item, $\partial(G)=V(G)$.
\item
If every vertex of $G$ is peripheral, then according to the previous item, $\kappa=n$.
Let $v$ a central vertex of $G$.
If every vertex of $G$ other than $v$ is peripheral, then $\partial(G)=V(G)-v$.
Let $v,w$ a pair of distinct central vertices of $G$. 
Then,  $v \in \partial(w)$ and $w \in \partial(v)$.
Hence, in this case, $\partial(G)=V(G)$.
\end{enumerate}
\vspace{-.4cm}
\end{proof}
\end{lemma}

\newpage
\begin{theorem}
If ${\rm diam}(G)=2$, then $G$ is a BDM graph.
\label{diam2}
\begin{proof}
According to Lemma \ref{b=n}, $\kappa=|\partial(G)|=n-1$ if and only if $G$ contains a unique central vertex.
Assume, w.l.o.g., that $\partial(G)=[n-1]$.
Let $G'$ a graph of order $n$ and boundary number $n-1$ such that $\hat{D}_{G'}=\hat{D}_G$.
In $G$ and $G'$ are not isomorphic, then  $N_G(n)=[n-1]$ and for some vertex $k\in [n-1]$, $kn\not\in E(G')$.
Thus, $d_{G'}(k,n)\ge2$, which means that $n\in \partial(k)$, a contradiction.
\end{proof}
\end{theorem}

The previous result is far from being true for every graph $G=([n],E)$ of diameter ${\rm diam}(G)=3$,  even in the case $\kappa(G)=n-1$. 
For example, if $G=([n],E)$ denotes any of the graphs displayed in Figures \ref{n7},\ref{bip8},\ref{dh8} and \ref{split10}, then  ${\rm diam}(G)=3$, $\kappa(G)=n-1$, and it is not  a BDM graph.

As for the case $\kappa(G)=n-1$, we have obtained a sufficient condition for a graph of this type to be a BDM graph, which is shown next.

\begin{definition}
Let $e=xy$ be an edge of a graph $G$.
Consider the graph $G'=G-e$.
The edge $e$ is called \emph{$x$-irrelevant} if, 
for every vertex $z\in N(x)$ , $d_{G'}(y,z)\le 2$.
\end{definition}

%
%
%
%
%
%

\begin{theorem}
Given a graph $G=([n],E)$ be a graph such that  $\partial(G)=[n-1]$.
Then, either  $G$ is a BDM graph or there is a vertex $h\in [n]$ 
such that $h\in N(n)$  and  $e=nh$ is $n$-irrelevant.
\label{new1}
\begin{proof}
Suppose that, for every vertex $h\in N(n)$, the edge $e=nh$ is not $n$-irrelevant.
W.l.o.g. we assume that $N(n)=[h]$ and $d_{G'}(h,1)\ge3$.
Thus, $\hat{D}_{G'}$  and $\hat{D}_{G}$ are different.

\noindent
Let $e=nk$, being $h+1\le k \le n-1$.
Consider the graph $G'=G+e$.
We distinguish two cases.

\noindent
{\bf (1)}
Assume that $d_{G}(k,1)\ge3$.
Then, $\hat{D}_{G'}$  and the distance matrix of $\partial(G)$ in $G'$ are different, since $d_{G'}(k,1)=2$. 

\noindent
{\bf  (2)}
Suppose that for every $i\in[h]$, $d_{G}(k,i)\le2$.
Then, $n\in \partial_{G}(k)$, a contradiction.
\end{proof}
\end{theorem}

\begin{figure}[ht]
\begin{center}
\includegraphics[width=0.9\textwidth]{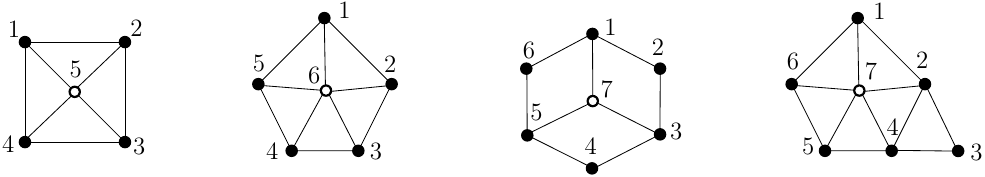}
\caption{
In all cases, the boundary is the set of black vertices,i.e., $\partial(G)=[n-1]$.}
\label{4irr}
\end{center}
\end{figure}

The converse of this result is far from being true. 
For example, for each of the graphs displayed in Figure \ref{4irr}, 
$\partial(G)=[n-1]$, $G$ is a BDM graph  and every edge incident to vertex $n$ is $n$-irrelevant.

\subsection{Ptolemaic graphs}\label{ptol}

A graph $G$ is called \emph{distance-hereditary} if every induced path is isometric.
A \emph{chordal graph} is a graph where there is no induced cycle of length at least four. 
A distance-hereditary graph is called \emph{Ptolemaic} if it is also  chordal.

\begin{theorem}[\cite{bm86}]
A graph $G$ is  Ptolemaic  if it can be constructed from a single vertex by the following operations:
\begin{enumerate}
\item[(1)] A leaf is added to an existing vertex.
\item[(2)] A vertex is duplicated to obtain two true twins.
\item[(3)] A vertex whose neighborhood is a clique is duplicated to get a false twin.
\end{enumerate}
\label{dhp}
\end{theorem}

We will use this construction to prove the following result.

\begin{lemma}\label{lem:ptolemaic}
The boundary $\partial(G)$ of a Ptolemaic graph $G$ consists of all its vertices except the cut-vertices.
\begin{proof}
We prove the above claim by induction in the sequence $G_1,G_2,\ldots ,G_n$ that leads to the graph $G$ according to the characterization shown in Theorem \ref{dhp}.
Obviously, the claim is true for $G_1=K_1$. 
Let us assume that it is also true for all $G_i$ with $1\leq i\leq k-1$ and we next prove it for $G_k$. 
For simplicity, $d_k(x,y)$ denotes the distance between vertices $x$ and $y$ in $G_k$, and $\partial_k(x)$ for the boundary of $x$ in $G_k$.

Suppose first that $G_k$ is obtained by adding a leaf $v$ to an existing vertex $u\in V(G_{k-1})$. 
Let $w,w'\in V(G_{k-1})$ such that $w,w'\neq u$. Since the distances from $w$ to the neighbors of $w'$ are the same in $G_{k-1}$ and in $G_k$ we have that $w'\in\partial_{k-1}(w)$ if and only if $w'\in\partial_k(w)$. Let $w\in V(G_{k-1})$ such that $w\in \partial_{k-1}(u)$. 
Now, it is clear that $w\in \partial_k(v)$. 
Moreover, if $w\notin \partial_{k-1}(u)$ and since the distances from $u$ to the neighbors of $w$ have not changed, then $w\notin \partial_k(u)$ and also $w\notin \partial_k(v)$ due to the fact that $d_k(v,x)=d_k(u,x)+1$ for any $x\in N[w]$. Finally, it is obvious that $v\in \partial_k(u)$. 
Hence, $\partial(G_k)$ contains all the vertices of $\partial(G_{k-1})$ except $u$, which is a cut-vertex,  plus $v$, which is not.

Suppose now, that $G_k$ contains one more vertex $v$ than $G_{k-1}$ which is a true twin of an existing vertex $u$. Note that $u\in \partial_k(v)$ and $v\in \partial(u)$, and that if $u$ is a cut-vertex in $G_{k-1}$, then it is not in $G_k$. The distances from $u$ to any other vertex $w\in V(G_{k-1})$ and between two vertices $w,w'\in V(G_{k-1})\setminus \{u\}$ have not changed so $w\in \partial(G_{k-1})$ if and only if $w\in\partial(G_k)$.

Finally, let us consider the case in which $G_k$ is obtained by adding a false twin $v$ to a vertex $u$ such that $N[u]$ is a clique. 
Then clearly, $u$ cannot be a cut-vertex in $G_{k-1}$ and $u\in\partial(G_{k-1})$, therefore $u$ is not a cut-vertex in $G_k$ and $u\in \partial(G_k)$ as well as $v$. 
Moreover, $u$ cannot be in a shortest path between two other vertices, so the distances between $w,w'\in V(G_{k-1})\setminus\{u\}$ have not changed and consequently $w\in \partial(G_{k-1})$ if and only if $w\in \partial(G_k)$. 
Hence, the proof is completed. 
\end{proof}
\end{lemma}

With this in mind, we are going to reconstruct a Ptolemaic graph $G$ from $\hat{D}(G)$ in two phases: first, we build its \emph{block-cutpoint tree} which is defined as the tree whose vertex set is the union of the set of blocks and the set of cutpoints of $G$, and in which two vertices are adjacent only if one corresponds to a block $B$ of $G$, and the other to a cutpoint $c$ of $G$ and $c|in B$ in $G$ (see~\cite{hp66}).

The second step would be then reconstruct all the blocks of the graph.

\begin{algorithm}
\caption{Reconstructing-Ptolemaic}
\label{algo1}
\begin{algorithmic}[1]
\Require A matrix $\hat{D}_{G}$ of a Ptolemaic graph $G$.
\Ensure The graph $G=(V,E)$.
\State Separate the vertices into 2-components; \label{step 1}
\State Compute the distances between any two 2-components $S$ and $S'$ as $d(S,S')=\min\{d(u,v): u\in S, v\in S'\}$;\label{step 2}
\State Construct a tree where any 2-component is a vertex and the distances between components are the same as in previous step; The resulted tree is the block-cutpoint tree of our graph;\label{step 3}
\State Rebuild the graph $G$.\label{step 4}
\end{algorithmic}
\end{algorithm}

\begin{corollary}
The Algorithm~\ref{algo1} runs in time $O(n^2)$.
\begin{proof}
Note that step~\ref{step 1} can be done with DFS in linear time. Also, the step~\ref{step 4} can be achieved in linear time by modifying the method of reconstructing a tree with the distances among the leaves given in~\cite{s62}. Finally, the step~\ref{step 4} can be done by adding one vertex at a time to the previous tree.

Clearly, the step that dominates the computation is step~\ref{step 2} which is done in time $O(n^2)$. 
\end{proof}
\end{corollary}

\begin{theorem}
Every Ptolemaic graph is a BDM graph
\begin{proof}
Let $G$ be a Ptolemaic graph with $n$ vertices and let $\hat{D}_G$ be its boundary distance matrix. 
According to Lemma~\ref{lem:ptolemaic}, that distance matrix contains the distance between all the vertices except the cut-vertices of the graph. 

First, we are going to divide the vertices of $\partial(G)$ in equivalence classes $S_1,\ldots ,S_p$ according to the following rules: $v\in S_i$ if and only if there exists $u\in S_i$ such that $d(v,u)=1$; if no such $u$ exists, then $v$ forms its own class of equivalence.

Clearly $G[S_i]$ is unique and consists of a block of $G$ without its cut-vertices. 
It only remains to construct the so called block-cutpoint tree of $G$. 
We can do that by applied the algorithm for trees given in~\cite{cp24} with a slight modification to consider non-leaves vertices. 
Once the tree is constructed, the vertices corresponding to each class $S_i$ are substituted by the vertices in $G[S_i]$ and the graph is reconstructed.
\end{proof}
\end{theorem}

\subsection{Interval graphs}\label{intg}

A graph is called an \emph{interval graph} if each of its vertices can be associated with an interval on the real line in such a way that two vertices are adjacent if and only if the associated intervals have a non-empty intersection. 

For any interval graph $G$, it is possible to give an enumeration $C_1,C_2,\ldots C_p$ of the cliques of $G$ such that if $x$ is a vertex such that $x\in C_i\cap C_k$, then $x\in C_j$ for all cliques with $i<j<k$ (see~\cite{fg65}).

\begin{figure}[htbp]
\begin{center}
\includegraphics[width=0.3\textwidth]{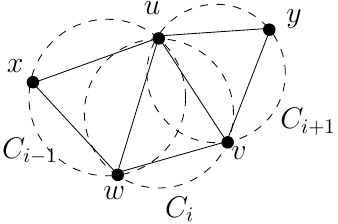}
\end{center}
\caption{How a 3-fan is constructed in an interval graph.}
\label{fig:interval1}
\end{figure}

With that characterization in mind, we prove the following:

\begin{theorem}
\label{lem:interval}
Let $G$ be an interval graph. 
Then, $u\in\partial(G)$ if and only if $u$ is neither a cut-vertex nor the central vertex of an induced 3-fan.
\end{theorem}
\begin{proof}
Before entering in the proof of this result, let us remember that if a vertex $x$ has a complete closed neighborhood, then $x\in\partial(G)$ for any graph $G$, not only interval graphs.

First, let $u\notin \partial(G)$. Since its closed neighborhood cannot be complete, then it should belong to at least two cliques, $C_i$ and $C_{i+1}$ (see Figure~\ref{fig:interval1}). 
Assume $u$ is not a cut-vertex, so there exists at least $v\in C_i\cap C_{i+1}$. 
If $N[u]\subseteq C_i\cup C_{i+1}$, then $u\in \partial(v)\subseteq \partial(G)$ contradicting our assumptions. 
Hence, there is a neighbor $w$ of $u$ belonging to $C_{i-1}$ or $C_{i+2}$ (suppose the first, without loss of generality) and $C_i$. Let $x\in C_{i-1}\setminus C_i$ different from $w$, and $y\in C_{i+1}\setminus C_i$. 
Then, $<u,v,w,x,y>$ induces a 3-fan having $u$ as central vertex, which leads to a contradiction.

Conversely, it is clear that if $u$ is a cut-vertex or the central vertex in a 3-fan, it does not belong to the boundary of $G$. 
\end{proof}

\begin{figure}[htbp]
\begin{center}
\includegraphics[width=0.3\textwidth]{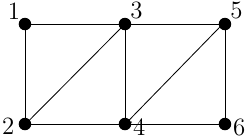}
\end{center}
\caption{Recognizing a 3-fan in an interval graph is not easy.}
\label{fig:interval2}
\end{figure}

Unlikely other graph families, having a characterization of the boundary of an interval graph does not lead immediately to be BDM-constructible (as an instance of the contrary, see Subsection~\ref{ptol}).  As an example of the difficulties of finding an algorithm to reconstruct an interval graph, consider the graph $G$ of Figure~\ref{fig:interval2}. It is an interval graph whose boundary consists of $\partial(G)=\{1,2,5,6\}$ as one can expect from Lemma~\ref{lem:interval}. However it is very difficult to detect the two 3-fans that are contained in the graph from the information given by $\hat{D}_G$.

\begin{conj}
\label{intervalconj}
Every interval graph is a BDM graph.
\end{conj}

\section{Some non-BDM graphs}\label{nbdm}

Through a brute-force algorithm applied to the list of all connected graphs of order at most $10$ (see~\cite{m}), we have obtained the following results.

\begin{figure}[ht]
\begin{center}
\includegraphics[width=0.9\textwidth]{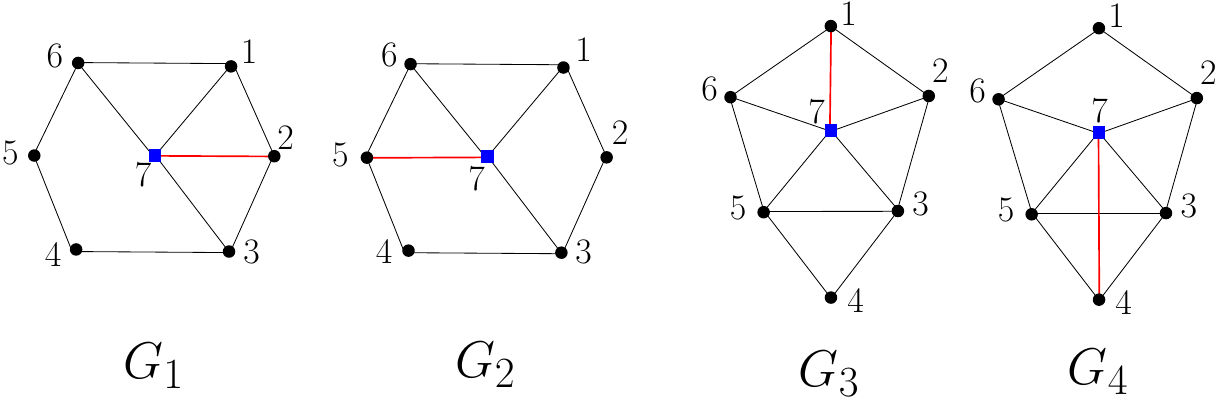}
\caption{These are the unique pairs of graphs of order $7$ having the same boundary distance matrix.
In all cases, $\partial(G_i)=[6]$.}
\label{n7}
\end{center}
\end{figure}

\begin{figure}[ht]
\begin{center}
\includegraphics[width=0.6\textwidth]{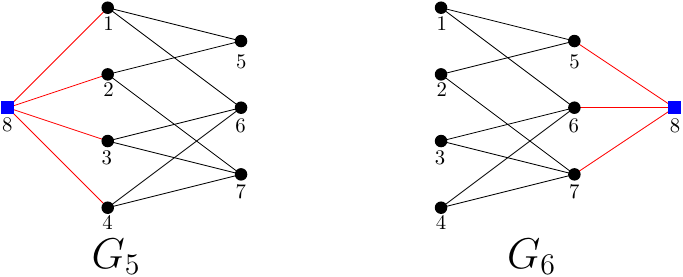}
\caption{This is the unique pair  of bipartite graphs of order $8$ having the same boundary distance matrix.
In both cases, $\partial(G_i)=[7]$.}
\label{bip8}
\end{center}
\end{figure}

\begin{figure}[ht]
\begin{center}
\includegraphics[width=0.6\textwidth]{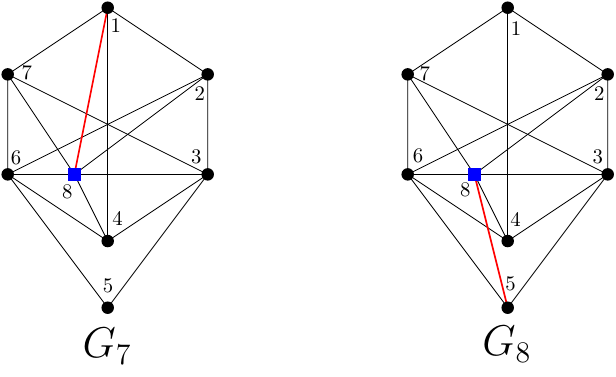}
\caption{This is the unique pair of distance-hereditary graphs of order $8$, having the same boundary distance matrix.
In both cases, $\hat{D}_{G_i}=[7]$.}
\label{dh8}
\end{center}
\end{figure}

\begin{proposition}
\label{mainconj7}
A   graph $G$ on $n$ vertices is a  BDM graph if at least one of the following conditions holds.
\begin{enumerate}[label=\rm \bf(\arabic*)]

\item
$n\le 7$ and $G\not\in \{G_1,G_2,G_3,G_4\}$ (See Figure \ref{n7}).

\item
$G$ is bipartite, $n\le 8$ and $G\not\in \{G_5,G_6\}$ (See Figure \ref{bip8}).

\item
$G$  is distance-hereditary, $n\le 8$ and $G\not\in \{G_7,G_8\}$ (See Figure \ref{dh8}).

\item
$G$  is chordal and $n\le9$.


\item
$G$  is an interval graph and $n\le10$.

\end{enumerate}

\end{proposition}

\begin{center} 
$\hat{D}_{G_1}=\hat{D}_{G_2}=$ 
\ma{rrrrrr} {0&1&2&3&2&1 \\ 1&0&1&2&3&2 \\ 2&1&0&1&2&2  \\ 3&2&1&0&1&2 \\ 2&3&2&1&0&1  \\ 1&2&2&2&1&0  }
\hspace{1.0cm}
$\hat{D}_{G_3}=\hat{D}_{G_4}=$ 
\ma{rrrrrr} {0&1&2&3&2&1 \\ 1&0&1&2&2&2 \\ 2&1&0&1&1&2  \\ 3&2&1&0&1&2 \\ 2&2&1&1&0&1  \\ 1&2&2&2&1&0  }
\end{center}

\begin{center} 
$\hat{D}_{G_5}=\hat{D}_{G_6}=$ 
\ma{rrrrrrr} 
{0&2&2&2&1&1&3 \\ 2&0&2&2&1&3&1 \\ 2&2&0&2&3&1&1  \\ 2&2&2&0&3&1&1 \\ 1&1&3&3&0&2&2  \\ 1&3&1&1&2&0&2  \\ 3&1&1&1&2&2&0 }
\end{center}

\begin{center} 
$\hat{D}_{G_7}=\hat{D}_{G_8}=$ 
\ma{rrrrrrrr} 
{0&1&2&1&3&2&1 \\
 1&0&1&2&2&1&2 \\
 2&1&0&1&1&2&1 \\
 1&2&1&0&2&1&2 \\
 3&2&1&2&0&1&2 \\
 2&1&2&1&1&0&1 \\ 
 1&2&1&2&2&1&0 }
\end{center}

\begin{figure}[ht]
\begin{center}
\includegraphics[width=0.55\textwidth]{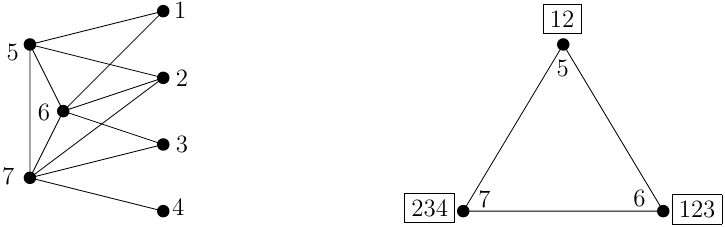}
\caption{A split graph of order $n=7$, independence number 4 and size $m=11$, drawn in two different ways.}
\label{split7}
\end{center}
\end{figure}

\begin{figure}[ht]
\begin{center}
\includegraphics[width=0.75\textwidth]{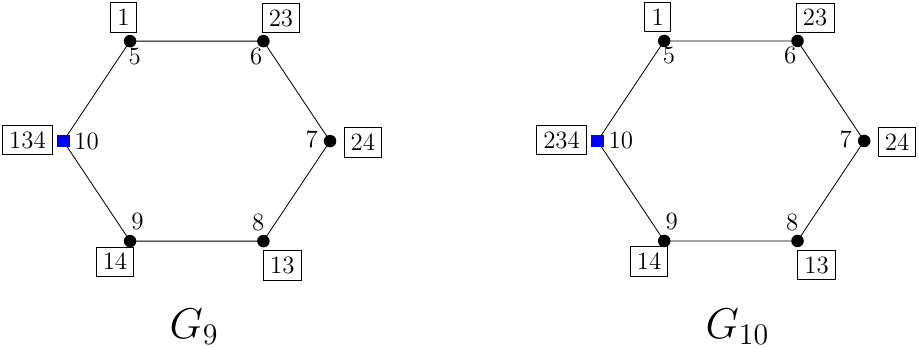}
\caption{A pair of split graphs of order $10$, having the same boundary distance matrix $\hat{D}_{G_i}=[9]$.}
\label{split10}
\end{center}
\end{figure}

A \emph{split graph} $G=(V,E)$ of order $n$ is a connected graph in which its vertex set $V$ can be partitioned into a clique $K_h$  and a (maximum) independent set $\bar{K}_{n-h}$.
It is denoted by $G=K_h\otimes \overline{K}_{n-h}$ (see Figure \ref{split7}).
Every split  graph $G$  has diameter at most diameter ${\rm diam}(G)=3$.
Certainly, as proved in Theorem \ref{diam2}.
Nevertheless, as shown in Figure \ref{split10}, not every split graph of diameter 3 is a BDM graph.

\begin{center} 
$\hat{D}_{G_9}=\hat{D}_{G_{10}}=$ 
\ma{rrrr|rrrrr} 
{0&3&2&2& 1&2&2&1&1 \\
 3&0&2&2& 2&1&1&2&2 \\
 2&2&0&2& 2&1&2&1&2 \\
 2&2&2&0& 2&2&1&2&1 \\
 \hline
 1&2&2&2& 0&1&1&1&1 \\
 2&1&1&2& 1&0&1&1&1 \\ 
 2&1&2&1& 1&1&0&1&1 \\ 
 1&2&1&2& 1&1&1&0&1 \\
 1&2&2&1& 1&1&1&1&0 }
\end{center}

\section*{Acknowledgement}
This work was supported by the Junta de Andalucía under Grant AGR-199; and UPC under Grant AGRUP-UPC.



\end{document}